\input amssym.def
\input amssym.tex
\magnification=1200
\font\titlefont=cmcsc10 at 12pt
\hyphenation{moduli}

%
%
\catcode`\@=11
\font\tenmsa=msam10
\font\sevenmsa=msam7
\font\fivemsa=msam5
\font\tenmsb=msbm10
\font\sevenmsb=msbm7
\font\fivemsb=msbm5
\newfam\msafam
\newfam\msbfam
\textfont\msafam=\tenmsa  \scriptfont\msafam=\sevenmsa
  \scriptscriptfont\msafam=\fivemsa
\textfont\msbfam=\tenmsb  \scriptfont\msbfam=\sevenmsb
  \scriptscriptfont\msbfam=\fivemsb
\def\hexnumber@#1{\ifcase#1 0\or1\or2\or3\or4\or5\or6\or7\or8\or9\or
      A\or B\or C\or D\or E\or F\fi }

\font\teneuf=eufm10
\font\seveneuf=eufm7
\font\fiveeuf=eufm5
\newfam\euffam
\textfont\euffam=\teneuf
\scriptfont\euffam=\seveneuf
\scriptscriptfont\euffam=\fiveeuf
\def\frak{\ifmmode\let\next\frak@\else
 \def\next{\errmessage{Use \string\frak\space only in math mode}}\fi\next}
\def\goth{\ifmmode\let\next\frak@\else
 \def\next{\errmessage{Use \string\goth\space only in math mode}}\fi\next}
\def\frak@#1{{\frak@@{#1}}}
\def\frak@@#1{\fam\euffam#1}

\edef\msa@{\hexnumber@\msafam}
\edef\msb@{\hexnumber@\msbfam}
\mathchardef\square="0\msa@03
\mathchardef\subsetneq="3\msb@28
\mathchardef\ltimes="2\msb@6E
\mathchardef\rtimes="2\msb@6F
\def\Bbb{\ifmmode\let\next\Bbb@\else
\def\next{\errmessage{Use \string\Bbb\space only in math mode}}\fi\next}
\def\Bbb@#1{{\Bbb@@{#1}}}
\def\Bbb@@#1{\fam\msbfam#1}
\catcode`\@=12
%
%

\def\HH{{\Bbb H}}

\def\PP{{\Bbb P}}
\def\C{{\bf C}}

\def\H{{\bf H}}
\def\P{{\bf P}}
\def\Z{{\bf Z}}
\def\ST{{\rm Sing}(\Theta)}
\font\eighteenbf=cmbx10 scaled\magstep3

\def\mapdown#1{\Big\downarrow\rlap{$\vcenter{\hbox{$\scriptstyle#1$}}$}}

\def\tto{\longrightarrow}
\def\llongrightarrow{\relbar\joinrel\relbar\joinrel\rightarrow}
\def\lllongrightarrow{\relbar\joinrel\relbar\joinrel\relbar\joinrel\rightarrow}
\def\vandaag{\number\day\space\ifcase\month\or
 januari\or februari\or  maart\or  april\or mei\or juni\or  juli\or
 augustus\or  september\or  oktober\or november\or  december\or\fi,
\number\year}
\magnification\magstep1
\vglue 2.0cm
\centerline{\eighteenbf The Moduli Space of Abelian Varieties }
\bigskip
\centerline{\eighteenbf  and the Singularities of the Theta Divisor}
\bigskip
\vskip 2pc
\centerline{\titlefont Ciro Ciliberto \& Gerard van der Geer}
\bigskip\bigskip
\centerline{\bf Introduction}
\bigskip
\noindent
The object of study here is the singular locus of the theta divisor $\Theta$ of
a principally polarized abelian variety $(X,\Theta)$. The special case
when $(X, \Theta)$ is the Jacobian of a curve $C$ shows that this is meaningful:
singularities of $\Theta$ are closely related to the existence of special linear
systems on the curve $C$ and for curves of genus $g\geq 4$ the divisor $\Theta$
is always singular.  But for the general principally polarized abelian variety
the theta divisor $\Theta$ is smooth. In their pioneering work [A-M]
Andreotti and Mayer introduced in the moduli space ${\cal A}_g$ of principally
polarized abelian varieties the loci $N_k$ of those principally polarized
abelian varieties for which $\Theta$ has a $k$-dimensional singular locus:
$$
N_k=\{ [(X, \Theta)] \in {\cal A}_g : \dim {\rm Sing}(\Theta) \geq k
\}\qquad k\geq 0.
$$
(Warning: we shall use a slightly different definition of the $N_k$
in this paper.)
They were motivated by the Schottky problem of characterizing Jacobian
varieties among all principally polarized abelian varieties. They showed
that $N_0$ is a divisor and that the image of the moduli space of curves under
the Torelli map is an  irreducible component of $N_{g-4}$ for $g\geq 4$.
In a beautiful paper [M2] Mumford calculated the
cohomology class of the divisor $N_0$ and used it to show that the moduli
space ${\cal A}_g$ is of general type for $g\geq 7$. Another reason for 
interest in the loci $N_k$ might be the  study of cycles on the moduli space
${\cal A}_g$.  The theory of automorphic forms seems to suggest that there
exist many algebraic cycles, but it seems very difficult to find them. The $N_k$
yield many interesting examples of cycles. Unfortunately, for $k\geq 1$ our
knowledge of the codimension of the $N_k$, let alone of the irreducible
components of $N_k$, is very limited.
\par
In this paper we give a new result on the codimension of the $N_k$ and
formulate a conjectural lower bound for the codimension of the $N_k$ in ${\cal A}_g$.
\smallskip
The first author would like thank A.\ Verra for useful discussions and the second author would like to thank Roy Smith and Robert Varley
for useful comments on a first draft of this paper.
\eject
\bigskip
\centerline{\bf 1. Review of Known Results}
\bigskip
\noindent
Let $(X,\Theta_X)$ be a principally polarized abelian variety of dimension
$g$ over $\C$. We shall assume that $g\geq 2$. Then $\Theta=\Theta_X$ is an
ample effective divisor with $h^0(\Theta)=1$ and $\Theta$ defines an
isomorphism 
$$
\lambda: X {\buildrel \approx \over \longrightarrow} \hat{X}, \qquad
x \mapsto [\Theta-\Theta_x],
$$
of $X$ with its dual abelian variety  $\hat{X}$; here $\Theta_x$ stands for
the translate of $\Theta$ over $x$.

We denote by ${\cal A}_g$ the moduli space of principally polarized abelian
varieties over $\C$. This is an orbifold of dimension $g(g+1)/2$.

The basic objects that we are interested in here are the loci
$$
N_k := N_{g,k} = \overline{ \{ [(X,\Theta_X)] \in {\cal A}_g(\C) : \dim {\rm
Sing}(\Theta_X) = k\} }\qquad \quad (0\leq k \leq g-2)
$$
where the overline means Zariski closure. Note that this definition is
slightly different from the definition of Andreotti-Mayer and the one used
by Mumford.  Andreotti and Mayer introduced these loci in 1967 in relation
with the Schottky problem. To explain the relation, we consider the moduli
space ${\cal M}_g$ of irreducible, smooth complete curves of genus $g$ over
$\C$. Then there is the Torelli map
$$
t: {\cal M}_g(\C) \to {\cal A}_g(\C)
$$
given by
$$
[C] \mapsto [(X={\rm Jac}(C), \Theta\subset {\rm Pic}^{g-1}(C))],
$$
where $\Theta$ is the divisor of effective line bundles of degree $g-1$ in
${\rm Pic}^{g-1}(C)$ and this defines  a divisor  up to
translations, again denoted by $\Theta$,  on ${\rm Jac}(C)$.  Riemann
showed (see [A-M])  that for a Jacobian of genus $g\geq 4$ the singular
locus
${\rm Sing}(\Theta)$ has dimension $g-4$ unless the curve $C$ is
hyperelliptic; in that case ${\rm Sing}(\Theta)$ has dimension $g-3$ and
there is even an explicit description of ${\rm Sing}(\Theta)$:
$$
{\rm Sing}(\Theta)= g_2^1 + W_{g-3}^0 \subset {\rm Pic}^{g-1}(C)
$$
with $W_{g-3}^0$ the locus of effective divisor classes (line bundles)
of degree $g-3$. We define
$$
J_g := \overline{t(M_g)}, \quad \hbox{\rm the Jacobian locus}.
$$
It is an irreducible closed algebraic subset of ${\cal A}_g$ of dimension $3g-3$. 
We also need $ H_g$, the hyperelliptic locus in ${\cal M}_g$, and we put
$$
{\cal H}_g = \overline{t(H_g)}, \quad \hbox{\rm the hyperelliptic locus
in } {\cal A}_g.
$$
It is irreducible and of dimension $2g-1$. Now Riemann's results on the
dimension of ${\rm Sing}(\Theta)$ for Jacobians imply:
$$
J_g \subseteq N_{g,g-4}\quad \hbox{\rm and}\quad
{\cal H}_g \subseteq  N_{g,g-3}.
$$
\par
\proclaim (1.1) Theorem. (Andreotti-Mayer) i) The Jacobian locus $J_g$ is an
irreducible component  of $N_{g,g-4}$. ii) The hyperelliptic locus ${\cal
H}_{g}$ is an irreducible component of 
$N_{g,g-3}$.
\par
The fact that $N_0$ is a divisor was first noticed in an unpublished
version of the Andreotti-Mayer paper. It was also proved by Beauville
([B2]) for $g=4$ and Mumford observed that Beauville's proof works in the
general case.

Mumford calculated in [M2] the class of $N_0$ in the Chow group 
$CH^1({\tilde{\cal A}}_g^{(1)})$ for the canonical partial compactification
${\tilde{\cal A}_g}^{(1)}$. If $\lambda=\lambda_1$ denotes the first Chern class
of the determinant of the Hodge bundle  (a line bundle the sections of which are
modular forms) and $\delta$ is  the class of the boundary then the result is
this:
\par
\proclaim (1.2) Theorem. (Mumford) The class of $N_0$  in
$CH^1( {\tilde{ \cal A}_g}^{(1)})$ is given by
$$
[N_0]= \big( {(g+1)!\over 2} + g!\big) \lambda - { (g+1)! \over 12} \,
\delta.
$$ 
\par
A variation of the notion of Jacobian varieties is given by {\sl Prym
varieties}:  take a double \'etale cover
$$
\pi: \tilde{C} \to C
$$
of smooth irreducible curves, where $\tilde{C}$ has genus $2g+1$ and $C$
has genus
$g+1$. The morphism $\pi$ induces a norm map 
${\rm Nm= \pi_*}: {\rm Pic}(\tilde{C}) \to {\rm Pic}(C)$. We now look
at the restriction to the degree $2g$ part
$$
\tilde{J} \cong {\rm Pic}^{2g}(\tilde{C}) \, {\buildrel {\rm Nm } \over
\longrightarrow } \, {\rm Pic}^{2g}(C) \cong {J} 
$$
and define the Prym variety $P$ as the connected component of zero in
 the kernel 
of ${\rm Nm}$:
$$
P= P(\tilde{C}/C)= \ker {\rm Nm} (\tilde{J} \to J)^0.
$$
As it turns out we can also write it as
$$
P \cong \{ L \in {\rm Pic}^{2g}(\tilde{C}): {\rm Nm}(L)= K_C, h^0(L) \equiv
0 (\bmod
\, 2)\}.
$$
and this is an abelian variety of dimension $g$. A principal 
polarization on $P$ is provided by the divisor of  effective line bundles
$$
\Xi = \{ L \in P :  {\rm Nm}(L)= K_C, h^0(L) \equiv 0 (\bmod \, 2),
h^0(L) > 0
\}.
$$
\par
Let $R{\cal M}_{g+1}$ be the moduli space of such double covers $\tilde{C}\to
C$. It is an orbifold of dimension $3g$ with a natural map $R{\cal M}_{g+1}\to 
{\cal M}_{g+1}$ defined by forgetting the cover $\tilde{C}$ of $C$.
The Torelli map has an analogue for this situation, the Prym-Torelli map:
$$
p: R{\cal M}_{g+1} \longrightarrow {\cal A}_g \qquad (\tilde{C} \to C) \mapsto
P(\tilde{C}/C).
$$
\par
\proclaim (1.3) Theorem. (Friedman-Smith [F-S], Donagi [Do 1,2]) The morphism $p$
is dominant for $g\leq 5$; it is birational to its image for $g\geq 6$, but not 
injective. 
\par
Note that the non-injectivity follows immediately from Mumford's
description of the Prym varieties of hyperelliptic curves, see [M1].\par
We define a locus in ${\cal A}_g$:
$$P_g := \overline{ p(R{\cal M}_{g+1})}, \quad \hbox{\rm the Prym locus in
${\cal A}_g$}
$$
This new locus is of dimension $3g$ for $g\geq 5$ and it contains the Jacobian
locus:
\par
\proclaim (1.4) Theorem. (Wirtinger, Beauville [B2]) The Prym locus contains the
Jacobian locus: $J_g \subseteq P_g$.
\par
The classical result of Riemann on the singular locus of $\Theta$ for
Jacobians has an analogue for Prym varieties. The singular points of $\Xi$ 
are of two types. If $L \in P$ then we have
$L \in {\rm Sing}(\Xi)$ if and only if
\item{i)} $h^0(L) \geq 4$, or
\item{ii)}  $L$ is of the form $\pi^*(E) +M$, where
$M\geq 0$ and $h^0(E)\geq 2$.
\par
The singularities of type i) are called {\sl stable} and those of type ii) are
called  {\sl exceptional}.  Welters and Debarre proved that
the singular locus of the divisor $\Xi$ has dimension $\geq g-6$,
see [W2,De3]. It
follows from their work and that of  Debarre ([D4]) that for a generic Prym
variety every singular point of $\Xi$ is stable and ${\rm Sing}(\Xi)$ is
irreducible of dimension $g-6$ for $g\geq 7$, reduced of dimension $0$
for $g=6$  and empty if $g\leq 5$. Mumford 
showed that if $\dim{\rm Sing}(\Xi)\geq g-4$ then ${\rm
Sing}(\Theta)$ has an exceptional component and the curve $C$ is either
hyperelliptic, trigonal, bi-elliptic, a plane quintic or a genus $5$ curve with
an even theta characteristic. By work of Debarre we know that if $C$ is not a
$4$-gonal curve then $\dim {\rm Sing}_{exc}(\Xi) \leq g-7$ for $g\geq 10$. He also
gives a beautiful description of exceptional singular locus of $4$-gonal curves.
\par
\proclaim (1.5) Theorem. (Debarre the  [D4]) The Prym locus $P_g$ is an irreducible
component of
$N_{g,g-6}$ for  $g\geq 7$.
\par
This shows once more that the components of $N_{g,k}$ give geometrically 
meaningful cycles on the  moduli space.
\bigskip
What do we know about the structure of the loci $N_{g,k}$ ? Let us start 
with $k=0$.
\par
\proclaim (1.6) Theorem. (Debarre [D5]) The divisor $N_{g,0}=N_0$ consists of
two  irreducible components
for $g\geq 4$:
$$
N_0= \theta_{\rm null} + 2 N_0^{\prime}.
$$
\par
Some explanation is in order here. The generic point of the
irreducible divisor $\theta_{\rm null}$  corresponds to a polarized
abelian variety $(X,\Theta)$ for which  $\Theta$ has one singularity, a
double point at a point of order $2$ of $X$, while the generic point of
the irreducible divisor $N_0^{\prime}$  corresponds to an abelian variety
$(X,\Theta)$ where $\Theta$ has two singularities. Mumford has shown in
[M] how the divisor $N_0$ can be defined scheme-theoretically so that it
comes with multiplicities. The divisor $\theta_{\rm null}$ is given as the zero
divisor of the modular form given by the product of the $2^{g-1}(2^g+1)$ 
even thetanulls $\theta[{\matrix{ \epsilon_1\cr \epsilon_2\cr}}](\tau,z)$.

\par
\noindent
{\bf (1.7) Example}. If $g=4$ the component $N_0^{\prime}$ is the Jacobian locus
$J_4$ as Beauville showed.  For $g=5$ the component $N_0^{\prime}$ can be
identified with the locus of intermediate Jacobians of double covers of $\PP^3$
ramified along a quartic surface with $5$ nodes, cf.\  [S-V1],[D5].
\smallskip
\par
Mumford showed that for $k\geq 1$ none of the $N_k$ have codimension $1$:
$$
{\rm codim}_{{\cal A}_g}N_k > 1 \quad{\rm if}\quad k \geq 1.
$$
At the other extreme we have $N_{g,g-2}$. We call a principally polarized
abelian variety {\sl decomposable} if it is  a product of 
(positive-dimensional) principally polarized abelian varieties.
The singular locus of
$\Theta$ for a decomposable  abelian variety has codimension $2$.
There are natural maps
$${\cal A}_{i_1} \times \ldots \times {\cal A}_{i_r} \to {\cal A}_g, \qquad
([X_{i_1}],\ldots, [X_{i_r}])
\mapsto [X_{i_1}\times\ldots \times  X_{i_r}].
$$
We denote the image by ${\cal A}_{i_1,\ldots,i_r}$.
Let 
$$\Pi_g= \bigcup_{1\leq i \leq g/2} {\cal A}_{i,g-i}
$$ 
be the locus of decomposable abelian varieties in ${\cal A}_g$.  This is a 
closed algebraic subset whose components have  codimension
$\geq g-1$ in ${\cal A}_g$.
\par
\proclaim (1.8) Proposition. $N_{g,g-2}= \Pi_g$ with $\Pi_g$  the locus of 
decomposable abelian varieties.
\par
This proposition is a corollary to a fundamental result by Koll\`ar
and a result of Ein-Lazarsfeld.
\par
\proclaim (1.9) Theorem. (Koll\`ar [K]) The pair $(X,\Theta)$ is log canonical. 
\par
Koll\`ar's result implies that 
$$
\Theta^{(r)}:= \{ x \in \Theta : {\rm mult}_x(\Theta) \geq r \}
$$
has codimension $\geq r$ in $X$.  Moreover, Ein and Lazarsfeld prove in
[E-L] 
$$
{\rm codim}_X\Theta^{(r)} =r \iff \hbox{$X$ is decomposable as product of $r$
p.p.a.v.}
$$
and this then implies the proposition, cf.\ also [S-V 3]. Ein and Lazarsfeld also
proved that if 
$\Theta$ is irreducible then $\Theta$ is normal and the singularities are
rational.
\par
\bigskip
\smallskip
\noindent
In the following table we collect what is known about components
of the $N_k$. We do not give multiplicities.
\smallskip
%
\font\tablefont=cmr8
\def\quad{\hskip 0.6em\relax}
\def\quod{\hskip 0.6em\relax}
\def\vhop{
    height2pt&\omit&&\omit&&\omit&&\omit&&\omit&&\omit&&\omit&&\omit&\cr}
\centerline{\bf Table for Low Genera}
$$
\vcenter{
\tablefont
\lineskip=1pt
\baselineskip=10pt
\lineskiplimit=0pt
\setbox\strutbox=\hbox{\vrule height .7\baselineskip
                                depth .3\baselineskip width0pt}%
\offinterlineskip
\hrule
\halign{&\vrule#&\strut\quod\hfil#\quad\cr
\vhop
&$g\backslash N_k$&&${\cal A}_g$&&$N_0$&&$N_1$&&$N_2$&&$N_3$&&$N_4$&&$N_5$&\cr
\vhop
\noalign{\hrule}
\vhop
&2&&$A_2=J_2$&&$\Pi_2$&&$\emptyset$&&$\emptyset$&&$\emptyset$&&
$\emptyset$&&$\emptyset$&\cr
\vhop
\noalign{\hrule} 
\vhop
&3&&$A_3=J_3$&&${\cal H}_3$&&
$\Pi_3$&&$\emptyset$&&$\emptyset$&&$\emptyset$&&$\emptyset$&\cr
\vhop
\noalign{\hrule}
\vhop
&4&&$A_4=P_4$&&$\theta_{\rm null}+J_4$&&${\cal H}_4$
&&$\Pi_4$&&$\emptyset$&&$\emptyset$&&$\emptyset$&\cr
\vhop
\noalign{\hrule}
\vhop
&5&&$A_5=P_5$&&$\theta_{\rm null}+
N_0^{\prime}$&&$J_5+A+B+C$&&${\cal H}_5$
&&$\Pi_5$&&$\emptyset$&&$\emptyset$&\cr 
\vhop
\noalign{\hrule}
\vhop
&6&&$A_6=?$&&$\theta_{\rm null}+
N_0^{\prime}
$&&$?$&&$J_6+?$&&${\cal H}_6+?$&&$\Pi_6$&&$\emptyset$&\cr
\vhop
\noalign{\hrule}
\vhop
&7&&$A_7=?$&&$\theta_{\rm null}+
N_0^{\prime}
$&&$P_7+?$&&$?$&&$J_7+?$&&${\cal H}_7+?$&&$\Pi_7$&\cr
\vhop
\noalign{\hrule}
\vhop
}
\hrule
}
$$
\par
\bigskip
It is known by work of Debarre that the three irreducible components
of $N_{5,1}$ different from $J_5$ have dimensions $10$, $9$ and $9$, cf.\
 [D1,Do3]. Moreover, in [D1] Debarre constructed components of $N_{g,g-4}$
(resp.
$N_{g,g-6}$) (resp. $N_{g,g-8}$) of codimension $g'(g-g')$ for $2\le g'\le
g/2$ and $g\ge 5$ (resp. for $3\le g'\le g/2$ and $g\ge 7$) (resp.\ for
$4\le g'\le g/2$ and $g\ge 9$) and  these thus are part of the question
marks in the table at positions $N_{5,1}$, $N_{6,2}$, $N_{7,3}$ and
$N_{7,1}$.
\bigskip
\noindent
\centerline{\bf 2.  Bounds and a Conjecture on the Codimension}
\bigskip
\noindent
As the review of the preceding section may show,  very little is known
about the components of the loci $N_{g,k}$. Apart from Mumford's
estimate that ${\rm codim}\, N_k \geq 2$ for $k\geq 1$ we know almost 
nothing about  the codimension of the $N_{g,k}$. Our new results give some
lower bounds for the codimension. Debarre proved  in an unpublished
manuscript independently that  ${\rm codim} N_{g,k} \geq k+1$ for $k\geq
1$.
\par
\proclaim (2.1) Theorem. Let $g\geq 4$. Then for $k$ with $1 \leq k \leq g-3$ 
we have ${\rm codim}\, N_{g,k} \geq k+2$. \par
\smallskip
\par
\proclaim (2.2) Theorem. Let $g\geq 5$. If $k$ satisfies  $g/3 < k
\leq g-3$ then ${\rm codim}\, N_{g,k} \geq k+3$.
\par
The first theorem is sharp for $k=1$ and $g=4,5$. However, we do not expect that this is an
accurate description of reality and believe that Theorems (2.1) and (2.2) are
never sharp for $k=1$ and $g\geq 6$, or for $k\geq 2$.  We conjecture the
following much stronger bound.
\par
\smallskip
\proclaim (2.3) Conjecture. If $1 \leq k \leq g-3$ and if $M$ is an
irreducible component of $N_{g,k}$  whose generic point corresponds to a
simple abelian variety then ${\rm codim}\, M
\geq { k+2 \choose 2}$. Moreover, equality holds if and only if
$g= k+3$ (resp.\ $g=k+4$) and then $M={\cal H}_g$ (resp.\ $M=J_g$).
\par
Note that by work of Beauville and Debarre ([B2,D1]) the conjecture is true for
$g=4$ and $g=5$.
\par
We now describe some corollaries of this. Let $\pi: {\cal X}_g\to {\cal A}_g$ be the
universal family of principally polarized abelian varieties. The reader should
view this as a stack, or replace ${\cal A}_g$ by a fine moduli space, e.g.\ the moduli
space of principally polarized abelian varieties with a level
$3$ structure. We can view ${\cal X}_g$ as the universal family of pairs $(X,
\Theta)$. In it we can consider the algebraic subset $S_g$ where the morphism
$\pi|\Theta $
is not smooth. If we write ${\cal A}_g$ as the orbifold $\H_g/{\rm
Sp}(2g,\Z)$ with $\H_g$ the upper half plane and ${\cal X}_g$ as the
orbifold
$$
\H_g \times \C^g / {\rm Sp}(2g,\Z) \ltimes \Z^{2g}
$$
then $\Theta$ is given in $\H_g \times \C^g$ by the vanishing of
Riemann's theta function  $\theta(\tau,z)=0$, with
$$
\theta(\tau,z)= \sum_{n \in \Z^g} e^{\pi i n^t \tau n + 2\pi i n^t z}
$$
and $S_g$ is defined in $\H_g \times \C^g$ by the $g+1$ equations
$$
\theta=0, \quad {\partial \theta \over \partial z_i} =0\qquad 
i=1,\ldots,g.
$$ 
Therefore, $S_g$ has codimension $\leq g+1$. Theorem (2.1) implies that the
codimension is equal to~$g+1$:
\par
\proclaim (2.4) Theorem. Every irreducible component of $S_g=
{\rm Sing}(\Theta)
\subset {\cal X}_g$ has codimension $g+1$ in ${\cal X}_g$, hence $S_g$
is locally a complete intersection. \par
\noindent
\par
\noindent
{\sl Proof.}   Take an irreducible component $S$ of $S_g$ and
let $N$ be its image under the natural map $\pi: {\cal X}_g \to
{\cal A}_g$. We first assume that $N$ is not contained in $\Pi_g$. Suppose now
 that  $N$ is contained in $N_{g,k}$ for some
$k\geq 1$, and we may assume by (1.8) that $ k \leq g-3$. Then by Theorem (2.1) the  codimension of $N$ in ${\cal A}_g$ is at
least
$k+2$. This  implies that  the codimension of $S$ in ${\cal X}_g$ is at
least $g+2$, which is impossible. Hence generically, the fibres of
$\pi_{|S}: S \to N$ are $0$-dimensional and $N$ must have codimension
$\leq 1$ in ${\cal A}_g$. So $S$ maps dominantly to a component of $N_{g,0}$, a
divisor and we get ${\rm codim} \geq g+1$. Finally, if $N$ is
contained in $\Pi_g$ we observe that $\Pi_g$ has codimension $g-1$
and the fibres have dimension $g-2$ leading also to ${\rm codim}\geq g+1$
and this concludes the proof.
\par
\smallskip
This Corollary of Theorem (2.1)  was obtained independently by  Debarre
in an unpublished note.
\par
\proclaim
(2.5) Corollary. $N_0$ is a divisor properly containing $\cup_{k \geq 1} N_k$.
\par
This raises the problem about the respective positions for higher $N_k$.
\par
\proclaim (2.6)  Problem. Is it true that $N_k$ properly contains 
$\cup_{i \geq k+1} N_i$ ?
\par
For the generic point of a component  $N_{0}$ we know  the singularities 
of $\Theta$. In general we know almost nothing about the nature of the
singular locus ${\rm Sing}(\Theta)$ of a generic point of a component of 
$N_{g,k}$. For a discussion of the case $N_1$ we refer to Section 8.

\bigskip
\centerline{\bf 3. Deformation Theory and the Heat Equation}
\bigskip
\noindent
In this section we explain Welters' interpretation of the Heat
Equation for the theta function, cf.\  [W1].  The Heat Equation is one of the tools
for obtaining our estimates on the codimension.
\par
Let $(X,\Theta)$ be a principally polarized abelian variety of dimension
$g$. We denote the invertible $O_X$-module $O_X(\Theta)$ associated to
$\Theta$ by $L$. The space
${\rm Def}(X)$ of linear infinitesimal deformations of the algebraic
variety
$X$ has a well-known cohomological interpretation:
$$
{\rm Def}(X) \cong H^1(X,T_X),
$$
where $T_X$ denotes the tangent sheaf of $X$. The space of linear
infinitesimal deformations of the pair $(X,\Theta)$ or equivalently of
the pair $(X,L)$, where we consider $\Theta$ or $L$ up to translations
on $X$, is given by
$$
{\rm Def}(X,L) \cong H^1(X,\Sigma_L),
$$
where $\Sigma_L$ is the sheaf of germs of differential operators of order
$\leq 1$ on $L$ (sums of functions and derivations). Given now a section
$s\in \Gamma(L)$ we obtain a complex
$$
0 \to \Sigma_L {\buildrel d^1s \over \longrightarrow} L \to 0, \quad 
\hbox{\rm with}\quad d^1s: D \mapsto D(s)
$$
on $X$ given by associating to a differential operator $D$ the section
$D(s)$ of $L$. The cohomological interpretation of the space of linear
infinitesimal deformations of the triple $(X,L,s)$ is  
$$
{\rm Def}(X,L,s) \cong {\HH}(d^1s),
$$
the hypercohomology of the complex. Explicitly, it can be given as
follows: if $(X_{\epsilon}, L_{\epsilon},s_{\epsilon})$ is an
infinitesimal deformation, then on a suitable open cover $U_j[\epsilon]$
of $X_{\epsilon}$ the section $s_{\epsilon}$ is given as $s_j+\epsilon
\sigma_j$ with $\sigma_j-\sigma_i= \eta_{ij}(s)$, where $\eta_{ij}(s)$
is a cocycle whose class in $H^1(X,\Sigma_L)$ determines the deformation
$(X_{\epsilon},L_{\epsilon})$. So we obtain a 1-cocycle $(\{ \sigma_i\},\{
\eta_{ij}\}) \in C^0(U,L)\oplus C^1(U,\Sigma_L)$ of the total
complex associated with
$$

\matrix
{&C^0(U,\Sigma_L) & \longrightarrow & C^1(U,\Sigma_L)
&\longrightarrow&...\cr
  &\mapdown{d^1s}      &         &\mapdown{-d^1s}&   \cr
 &C^0(U,L) & \longrightarrow & C^1(U,L) &\longrightarrow & ...\cr}
$$
and we thus have an element of ${\HH}^1(d^1s)$. 
\par
The central point is now the following:
\noindent
\proclaim (3.1) Claim. An element of $H^0(X, {\rm Sym}_2T_X)$ determines
canonically a linear infinitesimal deformation of $(X,L)$ and $(X,L,s)$.
\par
This follows from the first connecting homomorphism of
the exact sequence of hypercohomology of the short exact sequence of
complexes
$$

\matrix
{0 &\tto &\Sigma_L & \tto & \Sigma^{(2)}_L &\tto&{\rm
Sym}_2T_X&\tto & 0\cr
&  &\mapdown{d^1s}      &         &\mapdown{d^2s}& &\mapdown{}  \cr
0&\tto &L & {\buildrel{\rm id} \over\llongrightarrow } & L
&\lllongrightarrow & 0 & 
\tto &0,
\cr}\eqno(1)
$$
where $\Sigma^{(2)}_L$ stands for 
differential operators of order $\leq 2$ on $L$ and ${\rm Sym}_2T_X$
is the {\sl subspace} of elements fixed by the involution $(x_1,x_2)\mapsto 
(x_2,x_1)$ on $T_X\otimes T_X$. 
We thus have the connecting homomorphism of the upper exact sequence of (1)
$$
b: H^0(X, {\rm Sym}_2T_X) \tto H^1(X, \Sigma_L)
$$
and the connecting homomorphism of the short exact sequence of complexes (1)
$$
\beta: H^0(X, {\rm Sym}_2T_X) \tto \HH^1(d^1s)
$$
such that $b=f \cdot \beta$ with
$f: \HH^1(d^1s) \tto H^1(X,\Sigma_L)$ the forgetful
map, and  we find the morphisms 
$$ 
H^0(X,{\rm Sym}_2T_X) \tto \HH^1(d^1s)\tto H^1(X,\Sigma_L) \tto
H^1(X,T_X).\eqno(2)
$$
But it is well-known that we can identify $H^0(X, {\rm Sym}_2T_X)$ with
$H^1(X,\Sigma_L)$: 
we have $$
H^1(X,T_X) \cong H^1(X,O_X) \otimes
T_{X,0} \cong T_{\hat{X},0} \otimes T_{X,0}
$$
and using the polarization $\lambda: X {\buildrel \sim \over \tto}
\hat{X}$ we see that the subspace
corresponding to deformations preserving the polarizations is 
$$
{\rm Sym}_2T_{X,0} \subset T_{X,0} \otimes T_{X,0}\cong T_{\hat{X},0}
\otimes T_{X,0}.
$$ 
The composition $H^0(X, {\rm Sym}_2T_X)\to {\HH}^1(d^1s) \to
H^1(X, \Sigma_L)$ is therefore an isomorphism.
 The first spectral sequence for the hypercohomology
gives us an exact sequence
$$
H^0(\Sigma_L) \to  H^0(L){\buildrel \alpha \over \tto}  \HH^1(d^1s) \to H^1(\Sigma_L)
\to H^1(L),
$$
where $\alpha(t)= (X=X_{\epsilon},L= L_{\epsilon}, s+t\epsilon)$.
 So we get an exact sequence
$$
0 \to H^0(L)/ \C \cdot s \to \HH^1(d^1s) \to
H^1(\Sigma_L) \to 0
$$
which shows that for principally polarized abelian varieties
the forgetful map 
$$
f: {\HH}^1(d^1s) \to H^1(X,\Sigma_L)\eqno(3)
$$ 
is also an isomorphism: {\sl every deformation
$(X_{\epsilon}, L_{\epsilon})$ of $(X,L)$ canonically determines a deformation
$s_{\epsilon}$ of $s$.} This is Welters' interpretation of the classical Heat
Equations. If we represent
$X$ as a complex torus $X= \C^g/\Lambda$ with $\Lambda= \Z^g\tau+ \Z^g
$,
$\tau \in \H_g$, $z \in \C^g$ and $\Theta$ as before by
$$
\theta(\tau,z)= \sum_{n \in \Z^g} e^{\pi i n^t \tau n + 2\pi i n^t z}
$$
then it satisfies the relation
$$
2\pi i ( 1 + \delta_{ij}) {\partial \theta \over \partial \tau_{ij}} =
{\partial^2 \theta \over \partial z_i \partial z_j}
$$
where $\delta_{ij}$ denotes the Kronecker delta. These are the classical ``Heat
Equations" for Riemann's theta function.
\par
\bigskip
\noindent
\centerline{\bf 4.  Singularities of Theta and Quadrics}
\bigskip
\noindent
The tangent cone of a singular point $x$ of $\Theta$ with multiplicity $2$ 
defines after  projectivization and translation to the origin a
quadric $Q_x$ in $\P^{g-1}=\P(T_{X,0})$. Another description is obtained
as follows.  The singular points of
$\Theta\subset X$ are the points $x$ where the map $d^1s:\Sigma_L \to L$ of
Section 3 vanishes. Replace now in diagram (1) all sheaves by their fibres
at $x$ and denote the resulting maps by the suffix $(x)$. Then $(d^1s)(x)=0$ and
diagram (1) implies then that at such points $x$ the map $(d^2s)(x)$ factors
through
$$
({\rm Sym}_2T_X)_x \tto L_x. \eqno(4)
$$ 
This gives an element of  $L_x \otimes {\rm Sym}_2(T_X)_x^{\vee} \cong
H^0({\rm Sym}^2(\Omega_X^1))\otimes L_x$. We can view this as an
equation $q_x$ for  the projectivized tangent cone $Q_x$ of $\Theta$
at $x$ (if the multiplicity of the point is $2$; otherwise it is zero).
Note that if $x \in \ST$ 
then $\HH^1((d^1s)(x))$ can be identified with $L_x$ and
we can identify (4) with 
$$
H^0(X,{\rm Sym}_2T_X) \tto ({\rm Sym}_2T_X)_x \tto \HH^1((d^1s)(x))= L_x.
$$
The map $(d^2s)(x): {\rm Sym}_2T_X \to L_x$ sends an element $w$ to $0$
if and only if $q_x(w)=0$, i.e.\ if and only if the quadric $q_x$ and
the dual quadric $w$ are orthogonal.
\par
Suppose that we have an element
$w
\in  H^0(X, {\rm Sym}_2 T_X)$ determining by (3) an element of
$\HH^1(d^1s)$ with corresponding  deformation
$(X_{\epsilon},L_{\epsilon},s_{\epsilon})$ of $(X,L,s)$. This is given by a
cocycle $(\sigma_i, \eta_{ij})$ representing an element of $\HH^1(d^1s)$. With
respect to a suitable covering $\{ U_i \}$ of $X_{\epsilon}$ we can write the
section $s_{\epsilon}$ as
$$
s_i+\sigma_i \epsilon.
$$
Identifying $\HH^1((d^1s)(x))$ with $L_x$ we see that the corresponding
element of $L_x$ is given by $\sigma_i(x)$. Suppose that $x \in X$
deforms to $x_{\epsilon}$. The condition that
$s_{\epsilon}(x_{\epsilon})=0$ can be translated as follows:
$$
s_i(x) + (v_x s_i + \sigma_i(x))\epsilon =0,
$$
where $v_x$ is the tangent vector to $X$ at $x$ corresponding to
$x_{\epsilon}$. Since $s_i(x)=0$ and $v_x s_i=0$ because
$x$ is a singular point, the condition is $\sigma_i(x)=0$, i.e.\
$q_x(w)=0$. We thus see:
\smallskip
\par
\noindent
\proclaim (4.1) Lemma. Let $x$ be a quadratic singularity of $\Theta$.
 The infinitesimal deformations of $(X,\Theta)$
which keep $x$  on $\Theta$  are the 
deformations contained in $Q_x^{\bot}\subset {\rm Sym}_2(T_X)$.
In particular, the deformations that keep $x$ a singular point of
$\Theta$ are contained in $Q_x^{\bot}$. \par
\smallskip
\par
Let $R$ be an irreducible component of the locus ${\rm Sing}^{(2)}(\Theta)$ 
of quadratic singularities of $\Theta$, where we are assuming that  ${\rm
Sing}^{(2)}(\Theta)$ is not empty.  We now consider the map 
$$
\phi: R \to \P({\rm Sym}_2(T_X)_x^{\vee})= \P(T_{{\cal A}_g,[X]}^{\vee})
\qquad x \mapsto Q_x
$$
given by associating to $x\in R$ the quadric $Q_x\subset \P^{g-1}$. We
identified the space ${\rm Sym}_2(T_{X,0})$ with the tangent space
$T_{{\cal A}_g,[X]}$ to the moduli space ${\cal A}_g$ at $[(X,\Theta)]$.
\par
Since $\theta=0$ and all derivatives $\partial_j\theta$ vanish
on $R$ the partial derivatives $\partial_i\partial_j\theta$ are sections of
$O_R(\Theta)$:
\par
\noindent
\proclaim (4.2) Proposition. The map $\phi$ is given by sections of
$O_R(\Theta)$.
\par
Another way to interpret this is using the exact sequence
$$
0 \to {\cal T}_{\Theta} \tto {\cal T}_{X|\Theta} \tto {\cal N}_{\Theta,X} \tto
{\cal T}_{\Theta}^1 \to 0,
$$
where ${\cal T}_{\Theta}\cong {\cal O}_{\Theta}^g$ is the tangent sheaf, ${\cal
N}_{\Theta,X}\cong {\cal O}_{\Theta}(\Theta)$ is the normal sheaf, and ${\cal
T}_{\Theta}^1 \cong {\cal O}_{\ST}(\Theta)$ is the first higher tangent sheaf of
deformation theory and the middle arrow sends $\partial/\partial z_i$ to
$\partial\theta / \partial z_i$. The induced Kodaira-Spencer map is
$$
\delta: T_{{\cal A}_g,X} \tto H^0({\cal T}_{\Theta}^1 )=
H^0({\cal O}_{\ST}(\Theta))
$$
which maps $\partial/\partial \tau_{ij}$ to $\partial \theta /\partial
\tau_{ij}$,  cf.\ [S-V4].
 In this interpretation, for a singular point $x$ 
 the deformation $v\in T_{{\cal A}_g,X}$
keeps the point $x$ on $\Theta$ 
 if and only if $\delta(v)(x)=0$. If we assume
for simplicity that $\ST={\rm Sing}^{(2)}(\Theta)$ then the image of $\delta$
is a linear system on $\ST$ and we thus find a map
$$
\sigma: \ST \tto \PP(H^0({\cal T}_{\Theta}^1 )^{\vee}){\buildrel \delta^{\vee}
\over \tto } \PP(T_{{\cal A}_g,X}^{\vee}).
$$
The Heat Equations tell us that this can be identified with the map $\phi$
that associates to $x \in \ST$ the quadric defined by $\sum_{i,j}
({\partial^2\theta/ \partial z_i \partial z_j})\, x_ix_j$.

\smallskip
It might happen that all singularities of $\Theta$ are of higher order. In
order to deal with this case we extend the approach to the partial derivatives
of $s= \theta$ which are sections of $L$ when restricted to singular points of
$\Theta$. We define  
$$
R^{(j)}:=\{ x \in X : m_x(s)\geq j\},
$$
with $m_x$ the multiplicity at $x$, the set of points of multiplicity
$\geq j$ of $\Theta$; so $R^{(0)}=X$, $R^{(1)}=\Theta$, etc. Suppose
that $s$ is a non-zero-section of $L$. Then any partial derivative
$\eta= \partial_v\, s$ ($v \in {\rm Sym}^{(j)}(T_X)$) of weight
$j$ defines a section of $L|R^{(j)}$. 
\par
If $\eta$ is a partial derivative of $\theta$ then it satisfies again
the Heat Equation
$$
2\pi i (1+\delta_{ij }){ \partial \eta\over \partial \tau_{ij} }=
{\partial^2 \eta \over \partial z_i \partial z_j}.
$$
The algebraic interpretation is as follows. Given a partial
derivative $\eta$ of weight $j$ we apply the formalism of Section~3
to $\eta$ and find a map
$$
{\rm Sym}_2T_{X,0}\to H^0(R^{(j)}, {\rm Sym}_2T_{R^{(j)}})
\to \HH^1(d^1 \eta) \to H^1(R^{(j)},\Sigma_{R^{(j)}})  
$$
We claim that $\eta$ satisfies the heat equation: any linear
infinitesimal deformation of $(X,L)$ determines canonically a
deformation $\eta_{\epsilon}$ of $\eta$. This can be deduced in a way very
similar to the earlier case by extending Welters' analysis.
\par

\bigskip
\noindent
\centerline{\bf 5. The Tangent Space to $N_k$}
\bigskip
\noindent
Instead of looking at  a component of $N_0$ (with its reduced structure)
we may look at the space $\tilde{N}_0$ of triples
defined by
$$
S_g=\tilde{N}_0= \{ [(X,\Theta, x)] : x \in \ST \} \subset {\cal
X}_g,
$$
where ${\cal X}_g$ is the universal abelian variety over ${\cal
A}_g$. This has to be taken in the sense of stacks or one has to work
with level structures. We have a natural map
$\pi: \tilde N_0 \to N_0$. By Lemma (4.1)  the image under $d\pi$ of
the Zariski tangent space of $\tilde N_0$ at a point $[(X,x)]$ is contained in
the  space
$q_x^{\bot} \subset {\rm Sym}_2(T_{X,0})$.
\par
We shall call an abelian variety $X$ simple if it does not
contain abelian subvarieties $\neq X$ of positive dimension. 
The reason to consider simple abelian varieties is that we then
can use the non-degeneracy of the {\sl Gauss map}:
\bigskip
\proclaim (5.1) Theorem.  If $Z \subset X$ is a positive-dimensional smooth
subvariety of a simple abelian variety  then
the span of the tangent spaces to $Z$ translated to the origin is not
contained in a proper subspace of $T_{X,0}$. 
\par
Suppose that for $(X, \Theta)$ we have ${\rm
Sing}^{(2)}(\Theta)\neq \emptyset$ and that $N_k$ is smooth
at $[(X,\Theta)]$. The Zariski tangent space to $N_k$ (with its
reduced structure) at $[(X, \Theta)]$ is contained in the subspace of
${\rm Sym}_2(T_{X,0})$  orthogonal to the linear span in ${\rm
Sym}^2(\Omega^1_X)$ of the quadrics $q_x$ with $x \in  
R$ for some $k$-dimensional irreducible subvariety $R$ of ${\rm
Sing}(\Theta)$.
\par
By sending $x \in R$ to the quadric $Q_x$ we get a
natural map 
$$
R \tto \P(N_{N_k}),\qquad x \mapsto Q_x
$$
with $\P(N_{N_k})$ the projectivized conormal space to $N_k$.
  The image quadrics have rank $\leq g-k$ because of the following
lemma.
\smallskip
\proclaim (5.2) Lemma.  The Zariski tangent space to
 ${\rm Sing}^{(2)}(\Theta)$ at a point  $x$ equals ${\rm  
Sing}(Q_x)$. 
 \par
\noindent
{\sl Proof.} In local coordinates $z_1,\ldots, z_g$ a local equation
of $\Theta$ at a point $x$ is
$$
f= q_{x} + {\hbox{\rm higher order terms.}}
$$
By putting $q_{x}= \sum a_{ij}z_iz_j$ we get for $v=(v_1,\ldots,v_g)$
that $f(z+v)= \sum_{i,j} a_{ij}v_j z_i + \ldots$  and we see
 $v \in {\rm Sing}(Q_x)$, i.e.\ $ \sum_{i,j} a_{ij}v_j z_i=0$,  is
equivalent to
 $f(x+v)$ having  no linear term, i.e.\ $v \in T_{{\rm Sing}^{(2)}
(\Theta),x}$.
 \par
\smallskip
\noindent
\proclaim (5.3) Proposition. Let $X$ be a simple principally polarized abelian
variety and let $S$ be an open part of a    component  of ${\rm Sing}^{(2)}(\Theta)$ where
the rank of $Q_x$ is constant, say $g-d$.  Then the map
$S \to {\rm Gras}(d,g)$, $x \mapsto {\rm vertex}(Q_x)$ has finite fibres.\par
\smallskip
\noindent
{\sl Proof.} We first note by (5.2) that the tangent space at a point $x$ to
the reduced variety $S_{\rm red}$ is contained in the vertex of $Q_x$.
 If $F$ denotes a fibre of the map $x \mapsto Q_x$
then the tangent spaces to $F$ are contained in the subspace which
is the vertex of the constant  $Q_x$. 
The result then follows from (5.1).
\par
\noindent
\proclaim (5.4) Proposition. Let $X$ be simple and let $x$ be a quadratic
singularity of
$\Theta$ and a smooth point of ${\rm Sing}^{(2)}(\Theta)$. The general deformation
$w \in Q_x^{\bot}$ preserves only finitely many singularities of ${\rm
Sing}^{(2)}(\Theta)$. \par
\noindent
{\sl Proof.} The  deformations $w \in Q_x^{\bot}$ preserving 
$y\in {\rm Sing}(\Theta)$ are $Q_x^{\bot} \cap Q_y^{\bot}$. Hence all
deformations preserving $x$ preserve $y$ if and only if $Q_x=Q_y$. But the 
map $x \mapsto {\rm vertex}(Q_x)$ has finite fibres. 
 $\square$
\par
\smallskip \noindent
{\bf (5.5) Example}. Let $C$ be a curve of genus $g$ and $L\in \Theta\subset
{\rm Jac}(C)$ a quadratic singularity of $\Theta$. This means that the linear
system $|L|$ defined by $L$ is a $g_{g-1}^1$, i.e. has degree $g-1$ and
projective dimension $1$.\par
{\bf i)} $C$ is hyperelliptic. Then $L$ is of the form $g_2^1+D$, where $D$
is a divisor of degree $g-3$ on $C$. The quadric $Q_L$ is then the cone
projecting the canonical image of $C$ from the span of $D$. The image $\Sigma$ of
the map $\phi:\ST \to \PP^{{g+1 \choose 2}-1}$ can be identified with the
quadratic Veronese $V: \PP^{g-3} \tto \PP^{{g-1\choose 2}-1}\subset
\PP^{{g+1\choose 2}-1}$. So the normal space to $N_{g,g-3}$ at $[{\rm Jac}(C),
\Theta]$ is the subspace spanned by $\Sigma\cong \PP^{{g-1 \choose 2}-1}$. Since
${\rm codim}_{{\cal A}_g}{\cal H}_g = {{g+1\choose 2}- {g-1 \choose 2}}$ we see that
${\cal H}_g$ is a component of $N_{g,g-3}$.
\par
{\bf ii)} $C$ is not hyperelliptic. By a theorem of M.\ Green the space
spanned  by the quadrics $Q_x$ with $x \in {\rm Sing}^{(2)}(\Theta)$
is the space of quadrics containing the canonical curve. Since the space of
quadratic differentials on $C$ has dimension $3g-3$ it follows that the
normal space to $N_{g,g-4}$ at $[{\rm Jac}(C), \Theta]$ has dimension
$g(g+1)/2- (3g-3)$. We thus see that $J_g$ is a component of $N_{g,g-4}$.
\par
\bigskip
\centerline{\bf 6.  A Result on Pencils of Quadrics}
\bigskip
\noindent
One of the ingredients of our proofs is a classical result
of Corrado Segre on pencils of quadrics. Judging from the reactions
of experts this theorem seems to have been completely forgotten.
\par
\smallskip
\proclaim (6.1) Theorem. (C.\ Segre, 1883) Let $L$ be a linear
pencil of singular
quadrics of rank $\leq  n+1-r$ in $\PP^n$ with $n\geq 2$ whose generic member
has rank $n+1-r$ (i.e.\ the vertex $\cong \PP^{r-1}$). We assume
that the vertex is not constant in this pencil. Then the Zariski closure 
of the generic vertex in this pencil
$$
V_L= \big( \bigcup_{{\rm rk}(Q)= n+1-r} {\rm Vertex}(Q)\big)^{-}
$$
is a variety of dimension $r$ and degree $m-r+1$ in a projective
linear subspace $\PP^m \subset \PP^n$ with $m \leq (n+r-1)/2$ and
$r\leq (n+1)/3$. \par

If $L$ is a pencil of quadric cones whose generic member has rank $n$ in
$\PP^n$  and such that the vertex does not stay fixed then the
Zariski closure of the  
union of the vertices of the rank $n$ quadrics is contained
in the base locus of the family  
and it is a rational normal curve of degree $m$ contained in a
linear subspace $\PP^m \subset \PP^n$ with $m\leq n/2$.

For the proof we refer to Segre [S, p.\ 488-490]. It would be desirable to
have an extension of this theorem to higher dimensional linear families
of quadrics.

\bigskip
\centerline{\bf 7. Sketch of the Proof}
\bigskip
\noindent
We now sketch a proof of Theorem (2.1). Let $M$ be an irreducible
component of $N_k$. We choose a smooth point $\xi \in M$ corresponding
to a pair $(X, \Theta)$.  Note that we may assume that the abelian
variety $X$ is simple since the loci of non-simple abelian varieties have
codimension $\geq g-1$ in ${\cal A}_g$ and $g-1 \geq k+2$ by our assumption on $k$. We
choose a $k$-dimensional subvariety $R$ of ${\rm Sing}(\Theta_X)$
which deforms. For
simplicity we start with the case when the generic point $x$ of $R$ is a
quadratic singularity of $\Theta$.
\par
The construction of the preceding section yields a rational map
$$
\eqalign{ \phi : R \longrightarrow \PP(&N^{\vee}_{M/{\cal A}_g}), \quad x
\mapsto Q_x\cr
&{\|}\cr
& \PP^{\nu}\cr
} 
$$
Here $N^{\vee}_{M/{\cal A}_g}$ is the normal space of the component $M$ in
${\cal A}_g$ and we are assuming that the codimension of $M$ in ${\cal A}_g$ 
is $\nu +1$. 
We also have the {\sl Gauss map}
$$
\gamma: R_{\rm smooth} \longrightarrow {\rm Gras}(k,g), \quad x \mapsto {\rm
vertex }(Q_x)=\PP(T_{R,x})\eqno(5)
$$
which associates to a smooth point of $R$ its projectivized tangent
space, or equivalently the vertex of the quadric $Q_x$. The non-degeneracy
of the Gauss map $\gamma: R \to {\rm Gras}(k,g)$ implies
as in (5.3) that the fibres of
$\phi$ must be zero-dimensional, and this gives immediately $\nu \geq k$.
\par 
To go further we assume that $\nu =k$. Again, using the Gauss map we see that
$\phi$ maps dominantly to $\PP^{\nu}$. We consider a pencil $L$ of quadrics,
i.e. a $\PP^1 \subset \PP^{\nu}$. By (5.2) these quadrics have rank $\leq g-k$.
\smallskip
The theorem of Segre implies that the Gauss map $\gamma$ restricted to
$\phi^{-1}(L)$ is degenerate since the vertices lie in a linear subspace of
$\PP^{g-1}$, contradicting (5.1). 
\par
If the generic point of $R$ has higher multiplicity, say $r$, then we 
apply the preceding to the partial derivatives $\partial_v\theta$
with $|v|=r-2$ instead of to the section $\theta$. These satisfy the
Heat Equations and we proceed with these as with $\theta$ before. This completes
the sketch of proof of Theorem (2.1).
\par
For the proof of Theorem (2.2) we let $M$ as before be a component of
$N_k$, we pick a point $\xi \in M$ corresponding to $(X, \Theta)$ with $X$
simple and let  $R$ be a $k$-dimensional subvariety of ${\rm Sing}(\Theta)$
which deforms.
Assume then that ${\rm codim}(M) = k+2$ in ${\cal A}_g$. We get again a
rational map of $R$ to a projective space
$$
\phi: R \longrightarrow \P^{k+1}
$$
whose image is a hypersurface $\Sigma$. We have to distinguish two cases:
\item{i)} Not every quadric  corresponding to a point of $\PP^{k+1}$ is
singular.\par
\item{ii)} The general quadric corresponding to a point of $\PP^{k+1}$
is singular, say of of rank $g-r$.\par
First we treat the case i). Consider the discriminant locus $\Delta
\subset \PP^{k+1}$. This is a hypersurface of degree $g$ containing the image
$\Sigma$ of $R$ with multiplicity at least $k$:
$$
\Delta = k\Sigma + \Phi.
$$
\par
In order to be able to apply Segre's result we use the following well-known
Lemma: 
\proclaim (7.1) Lemma. A hypersurface of degree $\leq 2n-3$ in $\PP^n$
contains a line. \par
So if the degree of the hypersurface $\Sigma$ satisfies $\leq 2k-1$ we have
again a line and we can apply Segre's result. Note that since the degree of
$\Delta$ equals $g$ we have $\deg(\Sigma)
\leq g/k$, so that $g/k \leq 2k-1$ suffices and this follows from $g/3 <
k$. This rules out case i) if 
$\dim({\rm Sing}^{(2)}(\Theta))\geq k$. If the generic singularity of $\Theta$
has higher order we apply the procedure to the higher derivatives as before.
\par
To treat the remaining case ii), where all quadrics parametrized by $\PP^{k+1}$
are singular note that $r \leq k$ because  our quadrics generically 
 have rank $g-k$ by (5.2) and this should be less than $g-r$, the
generic rank of the whole family $\PP^{k+1}$.

 If $r=k$ then by Segre's result we get
$k\leq g/3$, contrary to our assumption. If $r<k$
we shall use a refinement of Segre's theorem which says that
the number of quadrics in a pencil of quadrics in $\PP^n$ where the rank
drops equals $n+r-2m-1$ in the notation of (6.1). Here one has
to count a quadric with  multiplicity $d$ if the rank drops by
$d$. In our case this yields
$$
(k-r) \deg \Sigma \leq (g-1)+r-2m -1 \leq g-2-r
$$
since $m\geq r$.
We get $\deg \Sigma \leq (g-2-r)/(k-r)\leq 2k-1$
and this assures us that
$\Sigma$ contains a line.
With this line we can apply Segre's theorem to get
a contradiction to the non-degeneracy of the Gauss map.
This concludes our sketch of proof.
\par
\smallskip
A closer analysis shows that we can draw stronger
conclusions from the proof. If
$[X,\Theta]$ is a  point of an irreducible component  of $N_1$
such that $X$ is simple and $\dim\ST =1$ then $[X,\Theta]$  admits a
linear deformation in codimension $3$ at most. If ${\rm codim}(N)=3$ and
$[X,\Theta] \in {\rm Sing}(N)$ then the singularities of $\Theta$ must ``get
worse."

\smallskip
This approach to getting estimates on the codimension of components of $N_k$
is by no means exhausted. For example, if ${\rm codim}(N)=k+3$ then the image
of a component $R$ of $\ST$ under the Gauss map is a codimension $2$ variety
$\Sigma$  in $\PP^{k+2}$. We now can use the variety spanned by the secant lines
instead of $\Sigma$ and apply Segre's result to that. We hope to return to this
point in the future (joint work with A.\ Verra)
\bigskip
\centerline{\bf 8. An Approach to the Conjecture for $N_1$}
\bigskip
\noindent
We now restrict to the case of $N_1$. Then the codimension is at least
$3$ and this is sharp: the case when the codimension is $3$ occurs. 
\smallskip
\par
\noindent
{\bf (8.1) Example}. Let $g=4$. We consider the hyperelliptic locus ${\cal
H}_4$. If $X={\rm Jac}(C)$, the Jacobian of a hyperelliptic curve of genus $4$,
then the singular locus ${\rm Sing}(\Theta)=g_2^1 + W_1^0$ is a copy of $C$ as
explained above. The class of this curve in the cohomology is $\Theta^3/3!$. By
associating to each point $x\in {\rm Sing}(\Theta)$ the vertex of the singular
quadric $Q_x\subset \PP^3$ we obtain the Gauss map $C \to \Gamma=\PP^1$ and the
image $\Gamma$ is the rational normal curve of degree $3$ in $\PP^3$. The
quadrics containing $\Gamma$ form a net $\PP^2$ of quadrics. In general for a
net of quadrics the curve of vertices is a curve of degree $6$ in $\PP^3$ and
the discriminant curve of singular quadrics in $\PP^2$ is a curve $\Delta$ of
degree $4$. But in our case the map $\ST \to \Gamma$ is of degree
$2$ to a rational curve and $\Delta$ is a conic with multiplicity $2$.
\par
\noindent
{\bf (8.2) Example}. Let $g=5$. We consider the Jacobian locus $J_5$. If
$X={\rm Jac}(C)$, the Jacobian of a curve of genus $5$, then the singular locus
${\rm Sing}(\Theta)= W_{4}^1$ is a smooth curve $\tilde{D}$ of genus $11$ and
class $\Theta^4/4!$ if $C$ is not trigonal and with no semi-canonical pencils.
The quotient of $\tilde{D}$ under  the involution $-1$ is a curve $D$ and
$\tilde{D} \to D$ is a double unramified cover. The Gauss map
$\tilde{D}$ is the Prym canonical map of $\tilde{D}$ to $\PP^4$. The map
$\phi: \tilde{D}\to \PP^2$  is a map of degree $2$ to a plane
quintic $\Delta$.
\par
\smallskip

Our Conjecture says that a component of $N_1$ is of codimension $3$ if
and only if $g=4$ (resp.\ $g=5$) and the component in question is
${\cal H}_4$ (resp.\ $J_5$). A tentative approach to proving this
might be the following.

Take an irreducible component $N$ of $N_1$ and assume it has codimension $3$ in
${\cal A}_g$. Let $[X,\Theta]$ be a general point of $N$ and $R$ a 
$1$-dimensional component of $\ST$.
\smallskip

{\bf Step i)} Try to prove that the generic point of an irreducible
component of $R$ is a double points of $\Theta$. If it is not, we should be able
to prove that ${\rm codim}(N)$ is higher than~$3$.
\smallskip
{\bf Step ii)}   Assume that the general quadric $Q_x$ of the span of $\Sigma$ is
smooth; otherwise use Segre's result. Let $\Delta$ be the discriminant locus of
degree $g$ in $\P^2$. The map $\phi: R \longrightarrow \Sigma \subset
\Delta$ has degree $\geq 2$ since it factors through $-1$. Prove that
the degree is $2$. This seems difficult. It would imply that $\Theta \cdot R
\leq 2g$.
\smallskip
{\bf Step iii)} We now assume that the class of  $R$ is a
multiple $m\alpha$ of the minimal class $\alpha= \Theta^{g-1}/(g-1)! \in
H^2(X,\Z)$. If it is not, then
${\rm End}(X)\neq \Z$ and this implies that ${\rm codim}(N) \geq g-1$. By
the preceding step we now find $m \leq 2$. 
\smallskip
{\bf Step iv)} Apply now the Matsusaka-Ran criterion or a result of
Welters. This implies that $X$ is a Jacobian or a Prym variety. The cases
$N\neq {\cal H}_4, J_5$ can then be ruled out by the results of Beauville.

\bigskip

\par
\centerline{\bf References}
\bigskip
\noindent
[A-M] Andreotti, A., Mayer, A.L.: On period relations for
abelian integrals on algebraic curves. {\sl Ann.\ Sc.\
Norm.\ Pisa \bf 3} (1967), p.\ 189--238.
\smallskip
\noindent
[B1] Beauville, A.: Vari\'et\'es de Prym et jacobiennes
interm\'ediaires. {\sl Ann.\ Sci.\ \'Ecole Norm.\ Sup.\ \bf (4) 10}
(1977), p.\  309-391.
\smallskip
\noindent
[B2] Beauville, A.:  Prym varieties and the Schottky problem. 
{\sl Invent.\ Math.\  \bf 41} (1977), p.\ 149--196.
\smallskip
\noindent
[B3] Beauville, A.: Sous-vari\'et\'es sp\'eciales des vari\'et\'es de 
Prym. {\sl Compositio Math.\ \bf 45} (1982), p.\ 357--383.
\smallskip
\noindent
[D1] Debarre, O.: Sur les vari\'et\'es ab\'eliennes dont le diviseur theta est
singulier en codimension $3$. {\sl  Duke Math.\ J.\ \bf 57 }(1988), p.\
221--273. 
\smallskip
\noindent
[D2] Debarre, O.:  Sur les vari\'et\'es de Prym des courbes t\'etragonales.
{\sl Ann.\ Sci.\ \'Ecole Norm.\ Sup.\ \bf (4) 21} (1988), p.\  545--559.
\smallskip
\noindent
[D3] Debarre,
O.:  Sur le probl\`eme de Torelli pour les vari\'et\'es de Prym. {\sl Amer.\ J.\
Math.\ \bf 111} (1989),p.\ 111--134.
\smallskip
\noindent
[D4] Debarre, O.: Vari\'et\'es de Prym et ensembles d'Andreotti et
Mayer.  {\sl Duke Math.\ J.\ \bf 60} (1990), p.\ 599--630.
\smallskip
\noindent
[D5]  Debarre, O.: Le lieu des vari\'et\'es ab\'eliennes dont 
le diviseur th\^eta est singulier a deux composantes.  
{\sl Ann.\ Sci.\ \'Ecole Norm.\ Sup.\ \bf 25} (1992), p.\  687--707. 
\smallskip 
 \noindent
[D6] Debarre, O.: Sur le th\'eor\`eme de Torelli pour les
solides doubles quartiques.
{\sl Compositio Math.\ \bf 73} (1990), p.\ 161-187.
\smallskip
\noindent
[Do1] Donagi R.: The tetragonal construction. {\sl Bull.\ A.M.S. \bf 4} (1981), p.\
181--185.
\smallskip
\noindent
[Do2] Donagi R.: The Schottky Problem. In {\sl Theory of Moduli}, p.\ 84--137.
Lecture Notes in Math.\ 1337, Springer Verlag (1988).
\smallskip
\noindent
[Do3] Donagi, R.: The fibers of the Prym map. In {\sl Curves, Jacobians, and
abelian varieties (Amherst, MA, 1990)}, p.\ 55--125, Contemp. Math., 136, Amer.\
Math.\ Soc., Providence, RI, 1992.
\smallskip
\noindent
[E-L] L.\ Ein, R.\ Lazarsfeld: Singularities of the theta divisor 
and the birational geomtry of irregular varieties. {\sl 
 J.\ Amer.\ Math.\ Soc.\ \bf 10} (1997), p.\ 243--258.
\smallskip
\noindent
[F-S] Friedman, R.\ , Smith, R.: The generic Torelli theorem for the Prym
map. {\sl Invent.\ Math.\ \bf 67 } (1982), p.\ 473--490.
\smallskip
\noindent
[K] J.\ Koll\`ar: Shafarevich maps and automorphic forms. Princeton
University Press, 1995.
\smallskip 
\noindent
[M1] Mumford, D.: Prym varieties. I. In {\sl  Contributions to analysis (a
collection of papers dedicated to Lipman Bers)}, p. 325--350. 
Academic Press, New York, 1974. 
\smallskip
\noindent
[M2] Mumford, D.: On the Kodaira Dimension of the Siegel Modular
Variety. In: {\sl Algebraic Geometry, Open Problems}. Proceedings
Ravello 1982. Lecture Notes in Math.\ 997 Eds.\ C.\ Ciliberto,
F.\ Ghione, F.\ Orecchia, p.\ 348-376.
\smallskip
\noindent
[S] Segre, C.: Ricerche sui fasci di coni quadrici in uno
spazio lineare qualunque. {\sl Atti della R.\ Accademia delle
Scienze di Torino  \bf XIX},(1883/4), p.\ 692-710. 
\smallskip
\noindent
[S-V 1] Smith, R., Varley, R.: Components of the locus of singular
theta divisors of genus $5$. In {\sl Algebraic geometry, Sitges } (Barcelona),
1983, 338--416, Lecture Notes in Math.\ 1124, Springer, Berlin-New York, 1985. 
\smallskip
\noindent
[S-V 2] Smith, R., Varley, R.: Singularity theory applied to
$\Theta$-divisors. In {\sl  Algebraic geometry (Chicago, IL, 1989)}, p.\ 238--257,
Lecture Notes in Math., 1479, Springer, Berlin, 1991. 
\smallskip
\noindent
[S-V 3] Smith, R., Varley, R.: Multiplicity $g$ points on theta divisors.
{\sl Duke Math.\ J.\ \bf 82} (1996), p.\ 319--326. 
\smallskip
\noindent
[S-V 4] Smith, R., Varley, R.: Deformations of theta divisors
and the rank $4$ quadrics problem. {\sl Comp.\ Math.\ \bf 76}
(1990), p.\  367--398.
\smallskip
\noindent
[W1] Welters, G.\ E.:
Polarized abelian varieties and the Heat Equation. {\sl Compositio Math.\ \bf 
49}, (1983), no. 2, 173--194. 
\smallskip
\noindent
[W2] Welters, G.: A theorem of Gieseker-Petri type for Prym
varieties. {\sl Ann.\ Scient.\ Ec.\ Norm.\ Sup.\ \bf 18} (1985), p.\
671--683.
\end